\documentclass{amsart}\usepackage{ifthen}\newboolean{amsart}\setboolean{amsart}{true}

\usepackage{amsmath,url}

\ifthenelse{\boolean{amsart}}{
\newtheorem{lemma}{Lemma}[section]
\newtheorem{theorem}[lemma]{Theorem}

\newtheorem{conjecture}[lemma]{Conjecture}
\newtheorem{problem}[lemma]{Problem}

}{}

\DeclareMathOperator{\gap}{gap}

\def\oG{\overline{G}}

\def\cE{{\mathcal E}}

\newcommand{\abssubkw}{
\begin{abstract}
Let $Q(n,\chi)$ denote the minimum clique size an $n$-vertex graph can have
 if its chromatic number is $\chi$.
Using Ramsey graphs we give an exact, albeit implicit, formula for
 the case $\chi\geq (n+3)/2$. 
\end{abstract}

\ifthenelse{\boolean{amsart}}{\subjclass[2010]{05C69, 05C35, 05D10}}{\subclass{05C69, 05C35, 05D10}}
\keywords{clique number, chromatic number, Ramsey graphs}}

\begin{document}

\title{Large chromatic number and Ramsey graphs}

\ifthenelse{\boolean{amsart}}{

\author{Csaba Bir\'o}
\address{Department of Mathematics, University of Louisville, Louisville, KY 40292, USA}
\email{csaba.biro@louisville.edu}

\author{Zolt\'an F\"uredi}
\address{Department of Mathematics, University of Illinois at Urbana-Champaign, Urbana, IL 61801, USA\\R\'enyi Institute of Mathematics of the Hungarian Academy of Sciences, Budapest, P. O. Box 127, Hungary-1364}
\email{z-furedi@illinois.edu}
\thanks{Research supported in part by the Hungarian National Science Foundation OTKA, and by the National Science Foundation under grant NFS DMS 09-01276.}

\author{Sogol Jahanbekam}
\address{Department of Mathematics, University of Illinois at Urbana-Champaign, Urbana, IL 61801, USA}
\email{jahanbe1@illinois.edu}

}{\author{Csaba Bir\'o
\and Zolt\'an F\"uredi\thanks{Research supported in part by the Hungarian National Science Foundation OTKA, and by the National Science Foundation under grant NFS DMS 09-01276.}
 \and Sogol Jahanbekam}
\institute{C. Bir\'o \at Department of Mathematics, University of Louisville, Louisville, KY 40292, USA\\ \email{csaba.biro@louisville.edu} \and
Z. F\"uredi \at Department of Mathematics, University of Illinois at Urbana-Champaign, Urbana, IL 61801, USA\\
R\'enyi Institute of Mathematics of the Hungarian Academy of Sciences, Budapest, P. O. Box 127, Hungary-1364\\\email{z-furedi@illinois.edu}\and
S. Jahanbekam \at Department of Mathematics, University of Illinois at Urbana-Champaign, Urbana, IL 61801, USA\\\email{jahanbe1@illinois.edu}
}}

\ifthenelse{\boolean{amsart}}{\abssubkw\maketitle}{\maketitle\abssubkw}

\section{Preliminaries}

The clique number, the chromatic number, and the independence number of a graph
 $G=(V,\cE)$ are denoted by $\omega(G)$, $\chi(G)$, and $\alpha(G)$, respectively.
Intuitively, large chromatic number must imply large cliques.
Define
\[
  Q(n,c) := \min \{\omega (G) : |V(G)| = n \text{ and } \chi(G) =
  c\}.\]

It is obvious that $Q(n,n)=n$, (the only graph to investigate is $K_n$),
 it is not difficult to show that $Q(n,n-1)=n-1$ $(n\geq 2)$
 (the complement of the graph should be a star)
and that $Q(n,n-2)\leq n-3$ for $n\geq 5$ (remove a five-cycle $C_5$ from $K_n$).
Bir\'o~\cite{Bir-11-u} determined $Q(n,n-k)$ for $k\leq 6$, whenever $n$ is sufficiently large,
$n> n_0(k)$.

\begin{equation*}
\begin{array}{rc*7{c@{\quad}}} 
k&=&0& 1& 2& 3& 4& 5& 6 \\
Q(n,n-k)&=& n& n-1& n-3& n-4& n-6& n-7& n-9 \\
\end{array}
\end{equation*}
Based on these values he was tempted to conjecture that
 if $n$ is large enough, then
$Q(n,n-k) = n-2k+\lceil k/2\rceil$.
He also showed $Q(n,n-k)\geq n-2k+3$ for $k\geq 5$ and $n$ is large enough.
Jahanbekam and West~\cite{Jah-Wes} observed that $Q(n,n-k)$ is
 at most the conjectured value whenever $n\geq 5k/2$ and they also asked if this
 threshold on $n$ is both sufficient and necessary for equality.
Their constructions is the complement of $\lfloor k/2\rfloor$ vertex disjoint $C_5$'s
 and a path $P_3$ if $k$ is odd.
The aim of this note is to give an exact formula for $Q(n,n-k)$ for $n\geq 2k+3$.
Our results established the above conjecture for $k \leq 12$ and $k=14$ but disproved it for any other value of $k$.

\section{Ramsey graphs}

The Ramsey number $R(3,\ell)$ is the minimum
integer $R$ such that every graph on $n\geq R$ vertices has either three
 independent vertices or a clique of size $\ell$.
It is well-known (Kim~\cite{Kim-95} and Ajtai, Koml\'os and
Szemer\'edi~\cite{Ajt-Kom-Sze-80}) that
 there are constants $\gamma_1, \gamma_2>0$ such that
\begin{equation}\label{eq:R3}
  \gamma_1 \frac{\ell^2}{\log \ell}< R(3,\ell) <  \gamma_2 \frac{\ell^2}{\log \ell}
  \end{equation}
hold for every $\ell\geq 3$.
The first few values are also known (see, e.g., the survey~\cite{Rad-94-11})
\begin{equation}\label{eq:R3tab}
\begin{array}{rc*{11}{c@{\quad}}}
\ell&=& 1& 2& 3& 4& 5& 6 &7    &8 &9 &10 &11\\
R(3,\ell)&=& 1& 3& 6& 9& 14& 18& 23  &28 & 36  & 40-43 & 46-51 \\
\end{array}
\end{equation}
To state our result we need an inverse of this function.
Let
\[
  \omega(x):=\min \{ \omega(G) : |V(G)|=x \text{ and } \alpha(G)\leq 2 \}.
\]
We have $\omega(x)=\omega$ for $R(3,\omega)\leq x < R(3,\omega+1)$.
For $k\geq 1$ define

\begin{equation*}
  q(k):= \min 
    \sum_{i=1}^{s} \left(  \omega(2k_i+1)-1\right)
  \end{equation*}
 where the minimum is taken over all positive integers $k_1, \dots, k_s$ with
 $k_1+\dots +k_s=k$, $s\geq 1$.  
Also define $q(0):=0$. From the tableaux (\ref{eq:R3tab}) one can easily calculate the
first few values of $q$
\begin{equation}\label{eq:qtab}
\begin{array}{rc*{12}{c@{\quad}}}
k & = & 0 & 1\text{--}2 & 3 & 4 & 5\text{--}6 & 7\text{--}8 & 9\text{--}10 & 11\text{--}13\\
\omega(2k+1)-1 & = & 0 & 1 & 2 & 3 & 3 & 4 & 5 & 6\\
q(k) & = & 0 & 1 & 2 & 2 & 3 & 4 & 5 & 6\\
\\
k & = & 14\text{--}17 & 18\text{--}19 & 20 & 21\text{--}22\\
\omega(2k+1)-1 & = & 7 & 8 & ? & 9\\
q(k) & = & 7 & 8 & ? & 9\\
\end{array}
\end{equation}
It also follows from (\ref{eq:R3}) that there exist $\gamma_1', \gamma_2'>0 $ such that
 for $k\geq 3$ we have
\[ \gamma_1' \sqrt{k \log k} < q(k)< \gamma_2' \sqrt{k \log k}. \]

\begin{theorem}\label{th:main}
For $n\geq 2k+3$
\begin{equation*}
  Q(n,n-k)= n-2k+q(k).
  \end{equation*}
 \end{theorem}

\section{The chromatic gap}

The \emph{chromatic gap} is defined as
\[
  \gap(n):= \max\{ \chi(G)-\omega(G): |V(G)|=n \}.
  \]
Gy\'arf\'as,  Seb\H o,  and Trotignon~\cite{Gya-Seb-Tro-11}
 showed that $\gap(n)= \lceil n/2\rceil  -\omega(n)$ for almost all $n$.
Our results are closely related, we use similar tools, but ours can be considered
 as a strengthening of theirs because, obviously,
  $\gap(n)= \max \{ c - Q(n,c)\}$.

\section{Graphs with independence number $2$}

The aim of this section is to prove that in the definition of $q(k)$, we may suppose that $s\leq 3$.

\begin{theorem} \label{th:FCJ}
Let $k$ be a positive integer.
Then there is  an integer $s$, $s\leq 3$,  such that
 $q(k)=\sum_{i=1}^s{\omega \left(2k_i+1\right)-1}$ where $k_1+\cdots+k_s=k$, $k_i$'s are positive integers.
\end{theorem}


\begin{conjecture}
The previous statement holds with $s=2$.
Even more, $q(k)=\omega(2k+1)-1$ for all $k$ except for $k=4$.
\end{conjecture}


\begin{lemma}[Xiaodong Xu, Zheng Xie, S. Radziszowski \cite{Xu-Xie-Rad-04}]\label{th:XZR}
Assume that $\omega_1 \geq \omega_2\geq 1$.
Then we have
\begin{equation}\label{eq:XZR}
\left(R\left(3,\omega_1+1\right)-1\right)+\left(R\left(3,\omega_2+1\right)-1\right)+\omega_2
   \leq R\left(3,\omega_1+\omega_2+1\right)-1.
\end{equation}
\end{lemma}

Since this is our main tool, for completeness, we include their construction.
\begin{proof}
Consider two vertex disjoint graphs $G_1$ and $G_2$ with $|V(G_i)|=R(3,\omega_i+1)-1$, $\alpha(G_i)\leq 2$ and $\omega(G_i)=\omega_i$ for $i=1,2$.
Let $R:=\{r_1,\ldots, r_{\omega_2} \}$ be a set disjoint to $V(G_1)$ and $V(G_2)$
 and suppose that $V_1:=\{v_1,\ldots, v_{\omega_2}\}$ and  $U_2:=\{u_1,\ldots, u_{\omega_2}\}$ are forming
 cliques in $G_1$ and $G_2$ respectively.
Define the graph $H$ with $V(H)=V(G_1)\cup V(G_2)\cup R$ as follows.
$H|V(G_i)=G_i$, $H|(R\cup V_1)$ and $H|(R\cup U_2)$ are complete graphs of sizes $2\omega_2$.
Connect every vertex in $v\in V\left(G_1\right)$ to every vertex
 in $u\in V\left(G_2\right)$ except if $v\in V_1$ and $u\in U_2$.
Finally, we have that  $r_i$ and $v_i$ have the same neighbors in $V(G_1)\setminus V_1$,
 similarly $N_H(r_i)\cap (V(G_2)\setminus U_2)=N_{G_2}(u_i)\setminus U_2$.
We obtain $|V(H)|=|V(G_1)|+|V(G_2)|+\omega_2$,  $\omega(H)=\omega_1+\omega_2$ and $\alpha(H)\leq 2$.
\end{proof}

\subsection*{Proof of Theorem~\ref{th:FCJ}}
Assume that $s$ is the minimum integer such that
 $q(k)=\sum_{1\leq i\leq s}$ ${\left(\omega(2k_i+1)-1\right)}$ and $k_1+\ldots+k_s=k$, $k_i\geq 1$.
Let $\omega_i:=\omega(2k_i+1)$.
By definition
\begin{equation}\label{eq:i}
  2k_i+1 \leq R(3,\omega_i+1)-1\quad \text{for all }i.
  \end{equation}
We may suppose that $s>1$, otherwise there is nothing to prove.

Our first observation is, that the multiset  $k_1, \dots, k_s$ is \emph{not reducible}.
This means that we cannot replace a set of $k_i$'s by their sum, i.e., for any subset
  $L \subset \{ 1,2,\dots ,s\}$, $2\leq |L|\leq s$ we have that
  \begin{equation*}
 \sum_{i\in L} \left(\omega(2k_i+1) -1\right)= \left( \sum_{i\in L} \omega_i\right)-|L| <
  \omega\left( 2\left(\sum_{i\in L} k_i\right)+ 1 \right)-1.
    \end{equation*}
This implies
\begin{equation*}
  R\left(3, \left(\sum_{i\in L}\omega_i\right)-|L|+ 2 \right)
   \leq  2\left( \sum_{i\in L} k_i\right)+1,
    \end{equation*}
which together with (\ref{eq:i}) give
\begin{equation}\label{eq:R}
  R\left(3, \left(\sum_{i\in L}\omega_i\right)-|L|+ 2 \right)
   \leq  \left( \sum_{i\in L} R(3, \omega_i+1)\right)-2|L|+1.
    \end{equation}

Suppose, on the contrary, that $s\geq 4$ and $\omega_1\geq \omega_2\geq \omega_3\geq \omega_4$.
We have $\omega_4\geq \omega(3)=2$.
Using the Erd\H os--Szekeres inequality $R(r,s)\leq R(r-1,s)+R(r,s-1)$ with $r=3$, we get
\begin{align}
R\left(3,\omega_2+1\right)-1&\leq R\left(3,\omega_2\right)+\omega_2, \label{eq:11}\\
R\left(3,\omega_3+1\right)-1&\leq R\left(3,\omega_3\right)+\omega_3, \label{eq:12}\\
R\left(3,\omega_4+1\right)-1&\leq R\left(3,\omega_4\right)+\omega_4. \label{eq:13}
\end{align}
Repeated applications of Lemma~\ref{th:XZR} gives
 \begin{align}
(R(3,\omega_2)-1)+(R(3,\omega_4)-1)+
  \omega_4-1 &\leq R(3,\omega_2+\omega_4-1)-1, \label{eq:14}\\
(R(3,\omega_1+1)-1)+(R(3,\omega_3)-1)+
      \omega_3-1 &\leq R(3, \omega_1+\omega_3)-1, \label{eq:15}\\
\begin{split}
(R(3,\omega_1+\omega_3)-1)+(R(3,\omega_2+\omega_4-1)-1)&+
 \omega_2+\omega_4-2\\
 \leq R(3,\omega_1&+\omega_2+\omega_3+\omega_4-2)-1. \label{eq:16}
\end{split}
\end{align}
Substitute $L:=\{ 1,2,3,4\}$  into (\ref{eq:R}) and
 add to the seven inequalities (\ref{eq:R})--(\ref{eq:16}).
We obtain
\begin{equation}\label{eq:17}
   \omega_4\leq 3.
    \end{equation}
If here equality holds then equality must hold in each of the 7 inequalities we added up.
However, (\ref{eq:13}) does not hold with equality for $\omega_4=3$ since
 $R(3,4)-1=9-1< 6+3= R(3,3)+3$.

From now on, we may suppose that $\omega_4=2$.
Substitute this to (\ref{eq:14}) and add  (\ref{eq:11}) to (\ref{eq:14}).
After rearrangement we get $2\leq \omega_2$.
So in the case $\omega_2\geq 4$ the sum of the right hand sides of
 (\ref{eq:R})--(\ref{eq:16}) exceeds the left hand sides by at least
 two, contradicting (\ref{eq:17}).
We obtain $\omega_2\leq 3$.

In the case of $\omega_2=3$, taking $L=\{ 2,4\}$, $\omega_2=3$, $\omega_4=2$ into (\ref{eq:R})
 we get the contradiction
 $R(3,5)=14 \leq R(3,3)+R(3,4)-2\times 2+1=12$.

The last case to consider is $\omega_4=\omega_3=\omega_2=2$. Substitute $L=\{2,3,4\}$ into (\ref{eq:R}) to obtain the final contradiction. \qed

\section{Construction}
In the case of $n\geq 2k+3$ the upper bound $n-2k+q(k)$ for $Q(n,n-k)$
 follows immediately from the definition of $q$ and Theorem~\ref{th:FCJ}.
Suppose that $k=\sum _{i=1}^{s} k_i$ such that
 $q(k)=\sum_{i=1}^{s} (\omega(2k_i+1)-1)$ where $k\geq s\geq 1$ and $s\leq 3$.
There is a graph $H_i$ with $2k_i+1$ vertices and with $\omega(H_i)=\omega(2k_i+1)$
 and with no three independent vertices.
Define the $n$-vertex graph $G$ in the following way. Consider the
 vertex disjoint union of $H_1$, \dots, $H_s$ and $n-\sum_{i=1}^s (2k_i+1)$ extra vertices. Then put an edge between any two vertices $u,v$, unless $u,v\in H_i$ for some $1\leq i\leq s$.
We have that $\omega(G)= n-\sum_{i=1}^s (2k_i+1) + \sum_{i=1}^s \omega(H_i)$, as stated.

Concerning the chromatic number of $G$, $\alpha(H_i)\leq 2$ implies that $\chi(H_i)\geq k_i+1$,
 hence
\[\chi(G)= n-\sum_{i=1}^s(2k_i+1)+ \sum_{i=1}^s \chi(H_i)\geq n-\sum_{i=1}^s k_i=n-k.
  \]
To obtain an example with chromatic number exactly
  $n-k$ delete edges arbitrarily from $G$ one by one until we obtain a subgraph $G'$
 with $\chi(G')=n-k$.
Since edge deletion does not increase the clique number we have
\[
  Q(n,n-k)\leq \omega(G')\leq \omega(G)= n-2k+q(k).
\]


\section{Lower bound by induction}
We will use induction on $k$, the cases $k=0,1$ are easy.
From now on, we suppose that $k\geq 2$.
The definition of $q$ immediately implies that
\begin{equation*}
   q(a)+ q(k-a)\geq q(k)
   \end{equation*}
for every integer $0\leq a\leq q$. In particular we have
\begin{equation}\label{eq:q-2}
   1+q(k-2)\geq q(k)
  \end{equation}
(for $k\geq 2$) and for all $0\leq x\leq q$
\begin{equation}\label{eq:q-x}
 x + q(k-x)\geq q(x)+q(k-x)\geq q(k).
  \end{equation}

Let $G$ be an $n$-vertex graph with $\chi(G)=n-k$, $n\geq 2k+3$, $k\geq 2$.
We  will show a lower bound for $\omega (G)$.
We distinguish two cases.

\subsection*{Case 1} $\alpha (G)\geq 3$.
Let $S\subset V(G)$ be a three-element independent set.
The chromatic number of the restricted graph  $G\setminus S:=G|(V\setminus S)$
 is at least $\chi(G)-1$.
We have that $|V(G\setminus S)|-\chi(G\setminus S)\leq k-2$.
Since $n-3\geq 2(k-2)+3$ we can use induction to $G\setminus S$.
\begin{multline*}
\omega(G)\geq\omega(G\setminus S)\geq Q(n-3, (n-3)-(k-2))\\
  \geq (n-3)-2(k-2)+q(k-2)= n-2k+(q(k-2)+1).
\end{multline*}
Then, (\ref{eq:q-2}) yields the desired lower bound.

\subsection*{Case 2} $\alpha(G)=2$.
Consider $\oG$, the complement of $G$.
The chromatic number of $G$ is $|V(G)|$ minus the matching number of $\oG$, $\nu(\oG)$.
So $\nu(\oG)=k$.
According to the Berge-Tutte formula, more exactly by the Edmonds-Gallai
 structure theorem (see, e.g.,~\cite{Lov-Plu-MT}) we have that
 there exists a partition of $V(\oG)= V_0\cup V_1 \cup \dots V_{a+b+c}\cup X$
such that
\begin{itemize}
\item $V_0$ is the set of isolated vertices of $\oG$,
\item $\oG\setminus (X\cup V_0)$ has $a+b+c\geq 1$ components, namely  $V_1, \dots, V_{a+b+c}$,
\item the sets $V_1, \dots, V_a$ are singletons, we define $k_h=0$ for $1\leq h \leq a$,
\item the sizes of $V_{a+1}, \dots, V_{a+b}$ are odd, $|V_i|=2k_i+1$, $k_i\geq 1$ for $a< i\leq a+b$,
\item the sizes of $V_{a+b+1}, \dots, V_{a+b+c}$ are even, $|V_j|= 2k_j$, $k_j\geq 1$ for $a+b< j\leq a+b+c$,
\item the matching numbers $\nu(\oG|V_i)=k_i$ for all $1\leq i\leq a+b+c$,
\item $0\leq |X|\leq a+b$, and finally
\item $k=\nu(\oG)=\left(\sum_{i=1}^{a+b+c} k_i\right) + |X|$.
\end{itemize}

We obtain
\[
\begin{split}
\omega(G)&\geq |V_0|+\sum_{i=1}^{a+b+c} \omega(G|V_i)\\
   &\geq \left( n-|X|-a-\sum_{i=a+1}^{a+b} (2k_i+1) -\sum_{j=a+b+1}^{a+b+c}
   2k_j\right)\\
   &\qquad + a + \sum_{i=a+1}^{a+b} \omega (2k_i+1) + \sum_{j=a+b+1}^{a+b+c} \omega(2k_j)\\
   &=n-|X| -2\left( \sum_{i=a+1}^{a+b+c} k_i \right) - b +
    \sum_{i=a+1}^{a+b} \omega (2k_i+1) + \sum_{j=a+b+1}^{a+b+c} \omega(2k_j) \\
   &=(n-2k) + |X| + \sum_{i=a+1}^{a+b}(\omega(2k_i+1)-1) + \sum_{j=a+b+1}^{a+b+c} \omega(2k_j) \\
   &\geq (n-2k) + |X| + \sum_{i=a+1}^{a+b+c}(\omega(2k_i+1)-1)\\
   &\geq (n-2k) + |X| + q(k-|X|)  \geq n-2k+q(k).
\end{split}
\]
In the last step we used (\ref{eq:q-x}), and in the previous one we used the obvious inequality
 $\omega(x)+1\geq \omega (x+1)$, which holds for every $x\geq 0$.
This completes the proof of the lower bound for $\omega(G)$, and also the proof of the Theorem. \qed

\section{Open problems and related questions}

The original motivation of this research was an analogue problem for partially
ordered sets (posets).

A \emph{realizer} is a set of linear extensions of the
poset $P$, such that their intersection (as relations) is $P$. The minimum
cardinality of a realizer is the dimension of the poset, a central notion in
poset theory. The ``standard example'' $S_n$ is the poset formed by considering
the $1$-element subsets and the $n-1$ element subsets of a set of $n$ elements,
ordered by inclusion. It is well known that $\dim(S_n)=n$, but there are posets
of arbitrarily large dimensions without including even $S_3$ as a subposet.

Hiraguchi \cite{Hir-55} proved that the dimension does not exceed half of the
number of elements of the poset. Bogart and Trotter \cite{Bog-Tro-73} showed
that for large $n$, the only $n$-dimensional poset on $2n$ points is $S_n$. But
what happens if the dimension is slightly less than half the number of
elements, is not known. We conjecture the following.

\begin{conjecture}
For every $t<1$, but sufficiently close to $1$
there is a $c>0$, if a poset has $2n$ points, and its dimension is at
least $tn$, then it contains a standard example of dimension $cn$.
\end{conjecture}

It is frequently noted that poset problems can be translated to graph theory
problems and vice versa by changing chromatic numbers of graphs to dimension of
posets, and cliques in graphs to standard examples in posets. Note that the
above conjecture would translate to the following statement: For every $t<1$,
but sufficiently close to $1$ there is a $c>0$ such that if a graph has $n$
points, and its chromatic number is at least $tn$, then it contains a clique of
$cn$ points. This graph version is trivial for all $t>1/2$.

The behavior of the sequence $\{Q(n,\lceil tn\rceil)\}_{n=1}^\infty$ is also
interesting in case $t\leq 1/2$. On one hand, our result implies that
$Q(n,\lceil tn\rceil)\leq Q(n,\lceil (n+3)/2\rceil)\leq 4+q(\lceil
n/2\rceil)=O(\sqrt{n\log n})$. On the other hand,
$\lim_{n\to\infty}Q(n,\lceil tn \rceil)\to\infty$: it follows from the
result that a graph with no cliques of size $t$ and $n$ vertices
has independence number $\Omega(n^{1/t})$, which is a straightforward
consequence of the classical Erd{\H o}s--Szekeres bound for
the Ramsey numbers.

\begin{problem}
It would be interesting to clarify more precisely what is the
exact relation between the (inverse of the) corresponding
Ramsey number and $Q(n,\lceil tn\rceil)$.
\end{problem}

A significant first step was done by Liu~\cite{Liu-12} for the case when
$1/t$ is a fixed integer and $n \to \infty$. 

We have reduced the determination of $Q(n,n-k)$ to the classical Ramsey number problem
 $R(3,\ell)$ whenever $n\geq 2k+3$.
It seems that with a bit of more care one can lower the bound of $n$ to $n\geq 2k+2$.
But below $2k$ one (probably) needs to use $R(4,\ell)$, too.

In the tableaux (\ref{eq:qtab}) one can find only a single case when $q(k)$ and
 $\omega(2k+1)-1$ differ from each other, namely $k=4$.
We \emph{conjecture} that this is the only case, i.e., $q(k)=\omega(2k+1)-1$ for $k\geq 5$.

We can also observe that
$q(k)=\lceil k/2\rceil$ holds for $k=0,1,2, \dots, 12$ and for $k=14$
 but it does not hold for any other value.

\ifthenelse{\boolean{amsart}}{\section{Acknowledgements}}{\begin{acknowledgements}}
The authors express their gratitude to Doug West for his
 valuable comments. 
\ifthenelse{\boolean{amsart}}{}{\end{acknowledgements}}

\bibliography{bib,extra}

\begin{thebibliography}{10}
\providecommand{\url}[1]{{#1}}
\providecommand{\urlprefix}{URL }
\expandafter\ifx\csname urlstyle\endcsname\relax
  \providecommand{\doi}[1]{DOI~\discretionary{}{}{}#1}\else
  \providecommand{\doi}{DOI~\discretionary{}{}{}\begingroup
  \urlstyle{rm}\Url}\fi

\bibitem{Ajt-Kom-Sze-80}
Ajtai, M., Koml{\'o}s, J., Szemer{\'e}di, E.: A note on {R}amsey numbers.
\newblock J. Combin. Theory Ser. A \textbf{29}(3), 354--360 (1980)

\bibitem{Bir-11-u}
Bir{\'o}, C.: Large cliques in graphs with high chromatic number (2011).
\newblock ArXiv:1107.2630

\bibitem{Bog-Tro-73}
Bogart, K.P., Trotter, W.T.: Maximal dimensional partially ordered sets. {II}.
  {C}haracterization of $2n$-element posets with dimension $n$.
\newblock Discrete Math. \textbf{5}, 33--43 (1973)

\bibitem{Gya-Seb-Tro-11}
Gy{\'a}rf{\'a}s, A., Seb{\H o}, A., Trotignon, N.: The chromatic gap and its
  extremes (2011).
\newblock ArXiv:1108.3444

\bibitem{Hir-55}
Hiraguchi, T.: On the dimension of orders.
\newblock Sci. Rep. Kanazawa Univ. \textbf{4}(1), 1--20 (1955)

\bibitem{Jah-Wes}
Jahanbekam, S., West, D.B.:
  \url{http://www.math.uiuc.edu/~west/regs/chromcliq.html}

\bibitem{Kim-95}
Kim, J.H.: The {R}amsey number {$R(3,t)$} has order of magnitude {$t^2/\log
  t$}.
\newblock Random Structures Algorithms \textbf{7}(3), 173--207 (1995)

\bibitem{Liu-12}
Liu, G.: Minimum clique number, chromatic number, and {R}amsey numbers.
\newblock Electron. J. Combin. \textbf{19}, Research Paper 55, 10 pp.
  (electronic) (2012).
\newblock
  \urlprefix\url{http://www.combinatorics.org/ojs/index.php/eljc/article/view/%
v19i1p55}

\bibitem{Lov-Plu-MT}
Lov{\'a}sz, L., Plummer, M.D.: Matching theory, \emph{North-Holland Mathematics
  Studies}, vol. 121.
\newblock North-Holland Publishing Co., Amsterdam (1986).
\newblock Annals of Discrete Mathematics, 29

\bibitem{Rad-94-11}
Radziszowski, S.P.: Small {R}amsey numbers.
\newblock Electron. J. Combin. \textbf{1}, Dynamic Survey 1, (electronic)
  (1994).
\newblock \urlprefix\url{http://www.combinatorics.org/Surveys/index.html}.
\newblock (version 2011 aug)

\bibitem{Xu-Xie-Rad-04}
Xu, X.D., Xie, Z., Radziszowski, S.P.: A constructive approach for the lower
  bounds on the {R}amsey numbers {$R(s,t)$}.
\newblock J. Graph Theory \textbf{47}(3), 231--239 (2004)

\end{thebibliography}
\ifthenelse{\boolean{amsart}}{\bibliographystyle{amsplain}}{\bibliographystyle{spmpsci}}

\end{document}